\documentclass[12pt,a4paper]{amsart}

\usepackage{amssymb}
\usepackage{graphicx}

\newcommand{\0}{\mathbf{0}}
\newcommand{\1}{\mathbf{1}}
\newcommand{\cle}{\preccurlyeq}
\newcommand{\rseth}{\widehat{\mathbf{R}}}
\newcommand{\rset}{\mathbf{R}}
\newcommand{\rmax}{\rset_{\max}}
\newcommand{\rmin}{\rset_{\min}}
\newcommand{\rmaxh}{\rseth_{\max}}
\newcommand{\smaxmin}{S_{\max,\min}}
\newcommand{\Mat}{\mathrm{Mat}}
\newcommand{\x}{\mathbf x}

\newcommand{\y}{\mathbf y}

\newcommand{\lx}{\underline{\x}}
\newcommand{\ly}{\underline{\y}}

\newcommand{\ux}{\overline{\x}}
\newcommand{\uy}{\overline{\y}}

\begin{document}

\title[Mathematics of semirings]{Universal algorithms, mathematics of
  semirings and parallel computations}

\author[G. Litvinov]{G. L. Litvinov${}^{*1,2}$}
\thanks{To appear in \textit{Springer Lecture
Notes in Computational Science and Engineering}.}
\thanks{${}^*$ Corresponding author: \texttt{glitvinov@gmail.com}}
\thanks{${}^1$ J.~V.~Poncelet Laboratory (UMI 2615 CNRS)}
\thanks{${}^2$  D.~V.~Skobeltsyn Research Institute for Nuclear Physics, Moscow State University}
\author[V. Maslov]{V. P. Maslov${}^3$}
\thanks{${}^3$ Physics Faculty, Moscow State University}
\author[A. Rodionov]{A. Ya. Rodionov${}^2$}
\author[A. Sobolevski]{A. N. Sobolevski${}^{1, 4}$}
\thanks{${}^4$ A.~A.~Kharkevich Institute for Information Transmission Problems}



\maketitle              
\index{Litvinov@G.L. Litvinov}
\index{Maslov@V.P. Maslov}
\index{Rodionov@A.Ya. Rodionov}
\index{Sobolevskii@A.N. Sobolevskii}

\begin{abstract}
This is a survey paper on applications of mathematics of semirings to numerical analysis and computing. Concepts of  universal algorithm and generic program are discussed. Relations between these concepts and mathematics of semirings are examined. A very brief introduction to mathematics of semirings (including idempotent and tropical mathematics) is presented. Concrete applications to optimization problems,  idempotent linear algebra and interval analysis are indicated. It is known that some nonlinear problems (and especially optimization problems) become linear over appropriate semirings with idempotent addition (the so-called idempotent superposition principle). This linearity over semirings is convenient for parallel computations.

\smallskip


\end{abstract}

\begin{center}
{\bf Contents}
\end{center}

\begin{itemize}
\item[] 1. Introduction
\item[] 2. Universal algorithms
\item[] 3. Universal algorithms and accuracy of computations
\item[] 4. Mathematics of semirings
	\begin{itemize}
		\item[] 4.1. Basic definitions
		\item[] 4.2. Closure operations
		\item[] 4.3. Matrices over semirings
            \item[] 4.4. Discrete stationary Bellman equations	
            \item[] 4.5. Weighted directed graphs and matrices over semirings
	\end{itemize}	
\item[] 5. Universal algorithms of linear algebra over semirings
\item[] 6. The idempotent correspondence principle
\item[] 7. The superposition principle and parallel computing
\item[] 8. The correspondence principle for computations
\item[] 9. The correspondence principle for hardware design
\item[] 10. The correspondence principle for software design
\item[] 11. Interval analysis in idempotent mathematics
\item[] References
\end{itemize}

\section{Introduction}

It is well known that programmers are the laziest persons in the world. They constantly try to write less and get more. The so-called art of programming is just the way to achieve this goal. They started from programming in computer codes, found this practice boring and invented assembler, then macro assembler, then programming languages to reduce the amount of work (but not the salary). The next step was to separate algorithms and data. The new principle ``Algorithm $+$ data structure $=$ program'' was a great step in the same direction. Now some people could work on algorithms and express those in a more or less data structure independent manner; others programmers implement algorithms for different target data structures using different languages.   

This scheme worked, and worked great, but had an important drawback: for the same algorithm one has to write a new program every time a new data type is required. Not only is this a boring process, it's also consuming time and money and worse, it is error-prone.  So the next step was to find how to write programs in a way independent of data structures. Object oriented languages, such as C++, Java and many others (see, e.g., [47, 62]) opened precisely such a way. For C++, templates and STL give even more opportunities.  This means that for a given algorithm one can vary a data structure while keeping the program unchanged.  The principle was changed to ``Algorithm $+$ data structure family $=$ program.'' 

But how does this approach work? What constitutes a ``data structure family''? The answer is this: operations. Any algorithm manipulates data by means of a set of basic operations. For example for sorting the comparison operation is required; for scalar products and matrix multiplications, the operations of addition and multiplication are required. So to use a new data structure with the same ``generic'' program one has to implement the required operations for this data.  

The hope was that some standard libraries of generic programs would cover almost all programmers' needs and programming would be very much like Lego constructing. Well, to some extent it worked, but not as much as it was expected.  Why? For different reasons. Some of them are economical. Why invest time and money in solving generic needs when a fast patch exists? Indeed, who would buy the next release if the current one is perfect?  Another reason: there is not a great variety of possible data structures for most popular algorithms. In linear algebra one can use floating numbers, double precision numbers, infinite precision floating numbers, rational numbers, rational numbers with fix precision, Hensel codes, complex numbers, integers. Not much.   

Mathematics of semirings gives a new approach to the generic programming. It parameterized algorithms. The new principle is ``Parameterized algorithm $+$ data structure family $=$ program.''  What are these parameters? They are operations (e.g., addition and multiplication). Sounds great, but does it really work? Yes. And we will show how and why. For example, we will show how the same old algorithm of linear algebra transformed by a simple redefinition of operations, can be applied to different problems in different areas. For example, the same program for solving systems of linear equations can be applied to the shortest path problem and other optimization problems (and interval versions of the problems) by means of simple redefinitions of the addition and multiplication operations. The algorithm was changed (by changing parameters, i.e. operations), the data structure was changed, the nature of the problem was changed, but the program was not changed!  

There are deep mathematical reasons (related to mathematics of semirings, idempotent and tropical mathematics) why this approach works. We will briefly discuss them in this paper. 


The concept of a generic program was introduced by many authors;
for example, in [36] such programs were called `program schemes.'
In this paper, we discuss {\it universal algorithms} implemented in the form of generic
programs and their specific features. This paper is closely related
to papers [38--40, 42, 43], in which the concept of a universal algorithm was
defined and software and hardware implementation of such algorithms was
discussed in connection with problems of idempotent mathematics
[31, 32, 37--46, 48--53, 55, 71, 72]. 
In the present paper the emphasis is placed on software and hardware implementations of universal algorithms, computation with arbitrary
accuracy, universal algorithms of linear algebra over semirings, and
their implementations.

We also present a very  brief introduction to mathematics of semirings and especially to
the mathematics  of \emph{idempotent semirings} (i.e.\ semirings with idempotent addition). Mathematics over idempotent semirings is called \emph{idempotent mathematics}. The so-called idempotent correspondence principle and idempotent superposition principle (see [38--43,  49--53] are discussed.

There exists a correspondence between interesting, useful, and
important constructions and results concerning the field of real (or
complex) numbers and similar constructions dealing with various idempotent
semirings. This correspondence can be formulated in the spirit of the
well-known N.~Bohr's \emph{correspondence principle} in quantum mechanics;
in fact, the two principles are intimately connected (see
[38--40]).  In a sense, the traditional
mathematics over numerical fields can be treated as a `quantum' theory, whereas the idempotent mathematics
can be treated as a `classical' shadow (or counterpart) of the traditional
one. It is important that the idempotent correspondence principle is valid for algorithms, computer programs and hardware units.

In quantum mechanics the \emph{superposition principle}  means that the
Schr\"odi\-n\-ger equation (which is basic for the theory) is linear.
Similarly in idempotent mathematics the (idempotent) superposition
principle (formulated by V.~P.~Maslov)
means that some important and basic problems and equations that are nonlinear in the usual sense (e.g.,
the Hamilton-Jacobi equation, which is basic for classical mechanics and appears in many 
optimization problems, or the Bellman equation and its versions and
generalizations) can be treated as
linear over appropriate idempotent semirings, see
[38--41, 48--51].

Note that numerical algorithms for infinite-dimensional linear problems
over idempotent semirings (e.g., idempotent integration, integral operators
and transformations, the Hamilton--Jacobi and generalized Bellman
equations) deal with the corresponding finite-dimensional approximations.
Thus idempotent linear algebra is the basis of the idempotent numerical
analysis and, in particular, the \emph{discrete optimization theory}.

B.~A.~Carr\'e [10, 11] (see also [7, 23--25]) used
the idempotent linear algebra to show that different optimization problems
for finite graphs can be formulated in a unified manner and reduced to
solving Bellman equations, i.e., systems of linear algebraic equations over
idempotent semirings.  He also generalized principal algorithms of
computational linear algebra to the idempotent case and showed that some of
these coincide with algorithms independently developed for solution of
optimization problems; for example, Bellman's method of
solving the shortest path problem corresponds to a version of Jacobi's
method for solving a system of linear equations, whereas Ford's algorithm
corresponds to a version of Gauss--Seidel's method. We briefly discuss Bellman equations and the corresponding optimization problems on graphs.

We stress that these well-known results can be interpreted as a
manifestation of the idempotent superposition principle.

We also briefly discuss interval analysis over idempotent and positive semirings. Idempotent internal analysis appears in [45, 46, 69]. Later many authors dealt with this subject, see, e.g. [12, 20, 27, 57, 58, 77].

It is important to observe that intervals over an idempotent semiring form a new idempotent semiring. Hence universal algorithms can be applied to elements of this new semiring and generate interval versions of the initial algorithms.

Note finally that idempotent mathematics is remarkably simpler than its traditional analog.

\section{Universal algorithms}

Computational algorithms are constructed
on the basis of certain primitive operations. These operations manipulate
data that describe ``numbers.'' These ``numbers'' are elements of a
``numerical domain,'' i.e., a mathematical object such as the field of
real numbers, the ring of integers, or an idempotent semiring of numbers
(idempotent semirings and their role in idempotent mathematics
are discussed in [10, 11, 22--26, 31, 32, 37--43] and below in this paper).

In practice elements of the numerical domains are replaced
by their computer representations, i.e., by elements of certain finite
models of these domains. Examples of models that can be conveniently used
for computer representation of real numbers are provided by various
modifications of floating point arithmetics, approximate arithmetics of
rational numbers [44], and interval arithmetics. The difference
between mathematical objects (``ideal'' numbers) and their finite
models (computer representations) results in computational (e.g.,
rounding) errors.

An algorithm is called {\it universal\/} if it is independent of a
particular numerical domain and/or its computer representation.
A typical example of a universal algorithm is the computation of the
scalar product $(x,y)$ of two vectors $x=(x_1,\dots,x_n)$ and
$y=(y_1,\dots,y_n)$ by the formula $(x,y)=x_1y_1+\dots+x_ny_n$.
This algorithm (formula) is independent of a particular
domain and its computer implementation, since the formula is
well-defined for any semiring. It is clear that one algorithm can be
more universal than another. For example, the simplest Newton--Cotes formula, the
rectangular rule, provides the most universal algorithm for
numerical integration; indeed, this formula is valid even for
idempotent integration (over any idempotent semiring, see below and [4, 32, 38--43, 48--51].
Other quadrature formulas (e.g., combined trapezoid rule or the Simpson
formula) are independent of computer arithmetics and can be
used (e.g., in the iterative form) for computations with
arbitrary accuracy. In contrast, algorithms based on
Gauss--Jacobi formulas are designed for fixed accuracy computations:
they include constants (coefficients and nodes of these formulas)
defined with fixed accuracy. Certainly, algorithms of this type can
be made more universal by including procedures for computing the
constants; however, this results in an unjustified complication of the
algorithms.

Computer algebra algorithms used in such systems as Mathematica,
Maple, REDUCE, and others are highly universal. Most of the standard
algorithms used in linear algebra can be rewritten in such a way
that they will be valid over any field and complete idempotent
semiring (including semirings of intervals; see below and [45, 46, 69], where
an interval version of the idempotent linear algebra and the
corresponding universal algorithms are discussed).

As a rule, iterative algorithms (beginning with the successive approximation
method) for solving differential equations (e.g., methods of
Euler, Euler--Cauchy, Runge--Kutta, Adams, a number of important
versions of the difference approximation method, and the like),
methods for calculating elementary and some special functions based on
the expansion in Taylor's series and continuous fractions
(Pad\'e approximations) and others are independent of the computer
representation of numbers.

\section{Universal algorithms and accuracy of computations}

Calculations on computers usually are based on a floating-point arithmetic
with a mantissa of a fixed length; i.e., computations are performed
with fixed accuracy. Broadly speaking, with this approach only
the relative rounding error is fixed, which can lead to a drastic
loss of accuracy and invalid results (e.g., when summing series and
subtracting close numbers). On the other hand, this approach provides
rather high speed of computations. Many important numerical algorithms
are designed to use floating-point arithmetic (with fixed accuracy)
and ensure the maximum computation speed. However, these algorithms
are not universal. The above mentioned Gauss--Jacobi quadrature formulas,
computation of elementary and special functions on the basis of the
best polynomial or rational approximations or Pad\'e--Chebyshev
approximations, and some others belong to this type. Such algorithms
use nontrivial constants specified with fixed accuracy.

Recently, problems of accuracy, reliability, and authenticity of
computations (including the effect of rounding errors) have gained
much atention; in part, this fact is related to the ever-increasing
performance of computer hardware. When errors in initial data and
rounding errors strongly affect the computation results, such as in ill-posed
problems, analysis of stability of solutions, etc., it is often useful
to perform computations with improved and variable accuracy. In
particular, the rational arithmetic, in which the rounding error is
specified by the user [44], can be used for this purpose.
This arithmetic is a useful complement to the interval analysis [54].
The corresponding computational algorithms must be
universal (in the sense that they must be independent of the computer
representation of numbers).

\section{Mathematics of semirings}

A broad class of universal algorithms is related to the concept of
a semiring. We recall here the definition of a semiring (see, e.g., [21, 22]).

 \subsection{Basic definitions}

Consider a semiring, i.e., a set~$S$ endowed with two
associative operations: \emph{addition}~$\oplus$ and
\emph{multiplication}~$\odot$ such that addition is commutative,
multiplication distributes over addition from either side, $\0$ (resp.,
$\1$) is the neutral element of addition (resp., multiplication), $\0 \odot
x = x \odot \0 = \0$ for all $x \in S$, and $\0 \neq \1$. Let the
semiring~$S$ be partially ordered by a relation~$\cle$ such that $\0$ is
the least element and the inequality $x \cle y$ implies that $x \oplus z
\cle y \oplus z$, $x \odot z \cle y \odot z$, and~$z \odot x \cle z \odot
y$ for all $x, y, z \in S$; in this case the semiring~$S$ is called
\emph{positive} (see, e.g., [22]).

A semiring~$S$ is called a \emph{semifield} if every nonzero element is invertible.

A semiring~$S$ is called \emph{idempotent} if $x \oplus x = x$ for all
$x \in S$, see, e.g., [4, 5, 7, 9, 10, 21--24, 31, 37--43, 48--53].
In this case the
addition~$\oplus$ defines a \emph{canonical
partial order}~$\cle$ on the semiring~$S$ by the rule: $x\cle y$ iff $x\oplus y=y$. It is easy to prove that any idempotent semiring is
positive with respect to this order. Note also that $x \oplus y =
\sup\{x,y\}$ with respect to the canonical order. In the sequel, we shall
assume that all idempotent semirings are ordered by the canonical partial
order relation.

We shall say that a positive (e.g., idempotent) semiring $S$ is \emph{complete} if it is complete  as an ordered set. This means that for every subset $T\subset S$ there exist elements $\sup T\in S$ and $\inf T\in S$.

The most well-known and important examples of positive semirings are
``numerical'' semirings consisting of (a subset of) real numbers and ordered
by the usual linear order $\leqslant$ on~$\rset$: the semiring~$\rset_+$
with the usual operations $\oplus = +$, $\odot = \cdot$ and neutral
elements $\0 = 0$, $\1 = 1$, the semiring~$\rmax = \rset \cup \{-\infty\}$
with the operations $\oplus = \max$, $\odot = +$ and neutral elements $\0 =
-\infty$, $\1 = 0$, the semiring $\rmaxh = \rmax \cup \{\infty\}$,
where $x \cle \infty$, $x \oplus \infty = \infty$ for all $x$, $x \odot
\infty = \infty \odot x = \infty$ if $x \neq \0$, and $\0 \odot \infty =
\infty \odot \0$, and the semiring~$\smaxmin^{[a,b]} = [a, b]$, where
$-\infty \leqslant a < b \leqslant +\infty$, with the operations $\oplus =
\max$, $\odot = \min$ and neutral elements $\0 = a$, $\1 = b$.  The
semirings~$\rmax$, $\rmaxh$, and~$\smaxmin^{[a,b]} = [a, b]$ are
idempotent. The semirings $\rmaxh$, $S^{[a,b]}_{\mathrm{max, min}}$, $\widehat{\textbf{R}}_+=\textbf{R}_+\bigcup\{\infty\}$ are complete.
Remind that every partially ordered set can be imbedded to its completion (a minimal complete set containing the initial one).

The semiring $\textbf{R}_\mathrm{min}=\textbf{R}\bigcup\{\infty\}$ with operations $\oplus=\mathrm{min}$ and $\odot=+$ and neutral elements $\textbf{0}=\infty$, $\textbf{1}=0$ is isomorphic to $\rmax$.

The semiring $\rmax$ is also called the \emph{max-plus algebra}. The semifields $\rmax$ and $\rmin$ are called \emph{tropical algebras}.
The term ``tropical'' initially appeared in [68] for a discrete version of the max-plus algebra as a suggestion of Ch.~Choffrut, see also [26, 55, 61, 72].

Many mathematical constructions, notions, and results over the fields of
real and complex numbers have nontrivial analogs over idempotent semirings.
Idempotent semirings have become recently the object of investigation of
new branches of mathematics, \emph{idempotent mathematics} and \emph{tropical geometry}, see, e.g. [5, 13, 15--19, 35--28, 31, 32, 37--46, 48--53, 55, 71, 72].

\subsection{Closure operations}

Let a positive semiring~$S$ be endowed with a partial unary \emph{closure
operation}~$*$ such that $x \cle y$ implies $x^* \cle y^*$ and $x^* = \1
\oplus (x^* \odot x) = \1 \oplus (x \odot x^*)$ on its domain of
definition. In particular, $\0^* = \1$ by definition. These axioms imply
that $x^* = \1 \oplus x \oplus x^2 \oplus \dots \oplus (x^* \odot x^n)$ if
$n \geqslant 1$. Thus $x^*$ can be considered as a `regularized sum' of the
series $x^* = \1 \oplus x \oplus x^2 \oplus \dots$; in an idempotent
semiring, by definition, $x^* = \sup \{ \1, x, x^2, \dots \}$ if this
supremum exists. So if $S$ is complete, then the closure operation is well-defined for every element $x\in S$.

In numerical semirings the operation~$*$ is defined as follows:
$x^* = (1-x)^{-1}$ if $x \prec 1$ in $\rset_+$, or $\widehat{\textbf{R}}_+$ and $x^*=\infty$ if $x\succcurlyeq 1$ in $\widehat{\textbf{R}}_+$; $x^* = \1$ if $x \cle \1$ in $\rmax$
and $\rmaxh$, $x^* = \infty$ if $x \succ \1$ in $\rmaxh$, $x^* = \1$
for all $x$ in $\smaxmin^{[a,b]}$. In all other cases $x^*$ is undefined.
Note that the closure operation is very easy to implement.

\subsection{Matrices over semirings}

Denote by $\Mat_{mn}(S)$ a set of all matrices $A = (a_{ij})$
with $m$~rows and $n$~columns whose coefficients belong to a semiring~$S$.
The sum $A \oplus B$ of matrices $A, B \in \Mat_{mn}(S)$ and the product
$AB$ of matrices $A \in \Mat_{lm}(S)$ and $B \in \Mat_{mn}(S)$ are defined
according to the usual rules of linear algebra:
$A\oplus B=(a_{ij} \oplus b_{ij})\in \mathrm{Mat}_{mn}(S)$ and
$$
AB=\left(\bigoplus_{k=1}^m a_{ij}\odot b_{kj}\right)\in\Mat_{ln}(S),
$$
where $A\in \Mat_{lm}(S)$ and $B\in\Mat_{mn}(S)$.
Note that we write $AB$ instead of $A\odot B$.

If the semiring~$S$ is
positive, then the set $\Mat_{mn}(S)$ is ordered by the relation $A =
(a_{ij}) \cle B = (b_{ij})$ iff $a_{ij} \cle b_{ij}$ in~$S$ for all $1
\leqslant i \leqslant m$, $1 \leqslant j \leqslant n$.

The matrix multiplication is consistent with the order~$\cle$ in the
following sense: if $A, A' \in \Mat_{lm}(S)$, $B, B' \in \Mat_{mn}(S)$ and
$A \cle A'$, $B \cle B'$, then $AB \cle A'B'$ in $\Mat_{ln}(S)$. The set
$\Mat_{nn}(S)$ of square $(n \times n)$ matrices over a [positive,
idempotent] semiring~$S$ forms a [positive, idempotent] semiring with a
zero element $O = (o_{ij})$, where $o_{ij} = \0$, $1 \leqslant i, j
\leqslant n$, and a unit element $E = (\delta_{ij})$, where $\delta_{ij} =
\1$ if $i = j$ and $\delta_{ij} = \0$ otherwise.

The set $\Mat_{nn}$ is an example of a noncommutative semiring if $n>1$.

The closure operation in matrix semirings over a positive semiring~$S$ can
be defined inductively (another way to do that see in [22] and below): $A^*
= (a_{11})^* = (a^*_{11})$ in $\Mat_{11}(S)$ and for any integer $n > 1$
and any matrix
$$
   A = \begin{pmatrix} A_{11}& A_{12}\\ A_{21}& A_{22} \end{pmatrix},
$$
where $A_{11} \in \Mat_{kk}(S)$, $A_{12} \in \Mat_{k\, n - k}(S)$,
$A_{21} \in \Mat_{n - k\, k}(S)$, $A_{22} \in \Mat_{n - k\, n - k}(S)$,
$1 \leqslant k \leqslant n$, by defintion,
\begin{equation}
\label{A_Star}
   A^* = \begin{pmatrix}
   A^*_{11} \oplus A^*_{11} A_{12} D^* A_{21} A^*_{11} &
   \quad A^*_{11} A_{12} D^* \\[2ex]
   D^* A_{21} A^*_{11} &
   D^*
   \end{pmatrix},
\end{equation}
where $D = A_{22} \oplus A_{21} A^*_{11} A_{12}$. It can be proved that
this definition of $A^*$ implies that the equality $A^* = A^*A \oplus E$ is
satisfied and thus $A^*$ is a `regularized sum' of the series $E \oplus A
\oplus A^2 \oplus \dots$.

Note that this recurrence
relation coincides with the formulas of escalator method of matrix
inversion in the traditional linear algebra over the field of real or
complex numbers, up to the algebraic operations used. Hence this algorithm
of matrix closure requires a polynomial number of operations in~$n$.

\subsection{Discrete stationary Bellman equations}

Let~$S$ be a positive semiring. The \emph{discrete stationary
Bellman equation} has the form
\begin{equation}
\label{AX}
	X = AX \oplus B,
\end{equation}
where $A \in \Mat_{nn}(S)$, $X, B \in \Mat_{ns}(S)$, and the matrix~$X$ is
unknown. Let $A^*$ be the closure of the matrix~$A$. It follows from the
identity $A^* = A^*A \oplus E$ that the matrix $A^*B$ satisfies this
equation; moreover, it can be proved that for positive semirings this
solution is the least in the set of solutions to equation (2)
with
respect to the partial order in $\Mat_{ns}(S)$.

\subsection{Weighted directed graphs and matrices over semirings}

Suppose that $S$ is a semiring with zero~$\0$ and unity~$\1$. It is well-known
that any square matrix $A = (a_{ij}) \in \Mat_{nn}(S)$ specifies a
\emph{weighted directed graph}. This geometrical construction includes
three kinds of objects: the set $X$ of $n$ elements $x_1, \dots, x_n$
called \emph{nodes}, the set $\Gamma$ of all ordered pairs $(x_i, x_j)$
such that $a_{ij} \neq \0$ called \emph{arcs}, and the mapping $A \colon
\Gamma \to S$ such that $A(x_i, x_j) = a_{ij}$. The elements $a_{ij}$ of
the semiring $S$ are called \emph{weights} of the arcs.

Conversely, any given weighted directed graph with $n$ nodes specifies a
unique matrix $A \in \Mat_{nn}(S)$.

This definition allows for some pairs of nodes to be disconnected if the
corresponding element of the matrix $A$ is $\0$ and for some channels to be
``loops'' with coincident ends if the matrix $A$ has nonzero diagonal
elements. This concept is convenient for analysis of parallel and
distributed computations and design of computing media and networks (see,
e.g., [4, 53, 73]).

Recall that a sequence of nodes of the form
$$
	p = (y_0, y_1, \dots, y_k)
$$
with $k \geqslant 0$ and $(y_i, y_{i + 1}) \in \Gamma$, $i = 0, \dots, k -
1$, is called a \emph{path} of length $k$ connecting $y_0$ with $y_k$.
Denote the set of all such paths by $P_k(y_0,y_k)$. The weight $A(p)$ of a
path $p \in P_k(y_0,y_k)$ is defined to be the product of weights of arcs
connecting consecutive nodes of the path:
$$
	A(p) = A(y_0,y_1) \odot \cdots \odot A(y_{k - 1},y_k).
$$
By definition, for a `path' $p \in P_0(x_i,x_j)$ of length $k = 0$ the
weight is $\1$ if $i = j$ and $\0$ otherwise.

For each matrix $A \in \Mat_{nn}(S)$ define $A^0 = E = (\delta_{ij})$
(where $\delta_{ij} = \1$ if $i = j$ and $\delta_{ij} = \0$ otherwise) and
$A^k = AA^{k - 1}$, $k \geqslant 1$.  Let $a^{(k)}_{ij}$ be the $(i,j)$th
element of the matrix $A^k$. It is easily checked that
$$
   a^{(k)}_{ij} =
   \bigoplus_{\substack{i_0 = i,\, i_k = j\\
	1 \leqslant i_1, \ldots, i_{k - 1} \leqslant n}}
	a_{i_0i_1} \odot \dots \odot a_{i_{k - 1}i_k}.
$$
Thus $a^{(k)}_{ij}$ is the supremum of the set of weights corresponding to
all paths of length $k$ connecting the node $x_{i_0} = x_i$ with $x_{i_k} =
x_j$.

Denote the elements of the matrix $A^*$ by $a^{(*)}_{ij}$, $i, j = 1,
\dots, n$; then
$$
	a^{(*)}_{ij}
	= \bigoplus_{0 \leqslant k < \infty}
	\bigoplus_{p \in P_k(x_i, x_j)} A(p).
$$

The closure matrix $A^*$ solves the well-known \emph{algebraic path
problem}, which is formulated as follows: for each pair $(x_i,x_j)$
calculate the supremum of weights of all paths (of arbitrary length)
connecting node $x_i$ with node $x_j$. The closure operation in matrix
semirings has been studied extensively (see, e.g., [1, 2, 5--7, 10, 11, 15--17, 22--26, 31, 32, 46]
and references therein).

\noindent\textbf{Example~1. The shortest path problem.}
Let $S = \rmin$, so the weights are real numbers. In this case
$$
	A(p) = A(y_0,y_1) + A(y_1,y_2) + \dots + A(y_{k - 1},y_k).
$$
If the element $a_{ij}$ specifies the length of the arc $(x_i,x_j)$ in some
metric, then $a^{(*)}_{ij}$ is the length of the shortest path connecting
$x_i$ with $x_j$.

\noindent\textbf{Example~2. The maximal path width problem.}
Let $S = \rset \cup \{\0,\1\}$ with $\oplus = \max$, $\odot = \min$. Then
$$
	a^{(*)}_{ij} =
	\max_{p \in \bigcup\limits_{k \geqslant 1} P_k(x_i,x_j)} A(p),
	\quad
	A(p) = \min (A(y_0,y_1), \dots, A(y_{k - 1},y_k)).
$$
If the element $a_{ij}$ specifies the ``width'' of the arc
$(x_i,x_j)$, then the width of a path $p$ is defined as the minimal
width of its constituting arcs and the element $a^{(*)}_{ij}$ gives the
supremum of possible widths of all paths connecting $x_i$ with $x_j$.

\noindent\textbf{Example~3. A simple dynamic programming problem.}
Let $S = \rmax$ and suppose $a_{ij}$ gives the \emph{profit} corresponding
to the transition from $x_i$ to $x_j$. Define the vector $B  = (b_i) \in
\Mat_{n1}(\rmax)$ whose element $b_i$ gives the \emph{terminal profit}
corresponding to exiting from the graph through the node $x_i$. Of course,
negative profits (or, rather, losses) are allowed. Let $m$ be the total
profit corresponding to a path $p \in P_k(x_i,x_j)$, i.e.
$$
	m = A(p) + b_j.
$$
Then it is easy to check that the supremum of profits that can be achieved
on paths of length $k$ beginning at the node $x_i$ is equal to $(A^kB)_i$
and the supremum of profits achievable without a restriction on the length
of a path equals $(A^*B)_i$.

\noindent\textbf{Example~4. The matrix inversion problem.}
Note that in the formulas of this section we are using distributivity of
the multiplication $\odot$ with respect to the addition $\oplus$ but do not
use the idempotency axiom. Thus the algebraic path problem can be posed for
a nonidempotent semiring $S$ as well (see, e.g., [66]). For
instance, if $S = \rset$, then
$$
	A^* = E + A + A^2 + \dotsb = (E - A)^{-1}.
$$
If $\|A\| > 1$ but the matrix $E - A$ is invertible, then this expression
defines a regularized sum of the divergent matrix power series
$\sum_{i \geqslant 0} A^i$.

There are many other important examples of problems (in different areas) related to algorithms of linear algebra over semirings (transitive closures of relations, accessible sets, critical paths, paths of greatest capacities, the most reliable paths, interval and other problems), see [1, 2, 4, 5, 10--13, 15--18, 20, 22--27, 31, 32, 45, 46, 53, 57, 58, 63--68, 74--77].

We emphasize that this connection between the matrix closure operation and
solution to the Bellman equation gives rise to a number of different
algorithms for numerical calculation of the closure matrix. All these
algorithms are adaptations of the well-known algorithms of the traditional
computational linear algebra, such as the Gauss--Jordan elimination, various
iterative and escalator schemes, etc. This is a special case of the idempotent superposition principle (see below).

In fact, the theory of the discrete stationary Bellman equation can be
developed using the identity $A^* = AA^* \oplus E$ as an additional axiom
without any substantive interpretation (the so-called \emph{closed
semirings}, see, e.g., [7, 22, 36, 66]).

\section{Universal algorithms of linear algebra over semirings}

The most important linear algebra problem is to solve the system
of linear equations
\begin{equation}
AX = B,
\end{equation}
where $A$ is a matrix with elements from the basic field and $X$ and
$B$ are vectors (or matrices) with elements from the same field.
It is required to find $X$ if $A$ and $B$ are given. If $A$ in (3)
is not the identity matrix $I$, then
system (3) can be written in form (2), i.e.,
$$
X = AX + B.\eqno{(2')}
$$
It is well known that the form (2) or ($2'$) is convenient for using the
successive approximation method. Applying this method with the initial
approximation $X_0=0$, we obtain the solution
\begin{equation}
X = A^*B,
\end{equation}
where
\begin{equation}
A^* = I+A+A^2+\cdots + A^n+\cdots
\end{equation}
On the other hand, it is clear that
\begin{equation}
A^* = (I-A)^{-1},
\end{equation}
if the matrix $I-A$ is invertible. The inverse matrix $(I-A)^{-1}$
can be considered as a regularized sum of the formal series (5).

The above considerations can be extended to a broad class of
semirings.

The closure operation for matrix semirings ${\Mat}_n(S)$ can be defined
and computed in terms of the closure operation for $S$ (see section 4.3 above); some
methods are described in [1, 2, 7, 10, 11, 22, 23--25, 32, 35, 42, 46, 65--67]. One such method is
described below ($LDM$-factorization).

If $S$ is a field, then, by definition, $x^*=(1-x)^{-1}$ for any $x\ne 1$. If $S$ is an idempotent semiring, then, by definition,
\begin{equation}
x^*=\1\oplus x \oplus x^2 \oplus\cdots=\sup\{\1, x, x^2, \dots\},
\end{equation}
if this supremum exists. Recall that it exists if $S$ is complete, see section~4.2.

Consider a nontrivial universal algorithm applicable to matrices over
semirings with the closure operation defined.

\begin{center}
{\it Example 5: Semiring $LDM$-Factorization.}
\end{center}

Factorization of a matrix into the product $A = LDM$, where $L$ and $M$
are lower and upper triangular matrices with a unit diagonal,
respectively, and $D$ is a diagonal matrix, is used for solving
matrix equations $AX = B$. We construct a similar
decomposition for the Bellman equation $X = AX \oplus B$.

For the case $AX = B$, the decomposition $A = LDM$ induces the following
decomposition of the initial equation:
\begin{equation}
   LZ = B, \qquad DY = Z, \qquad MX = Y.
\end{equation}
Hence, we have
\begin{equation}
   A^{-1} = M^{-1}D^{-1}L^{-1},
\label{AULinv}
\end{equation}
if $A$ is invertible. In essence, it is sufficient to find the
matrices $L$, $D$ and $M$, since the linear system (8) is easily
solved by a combination of the forward substitution for $Z$, the
trivial inversion of a diagonal matrix for $Y$, and the back
substitution for $X$.

Using (8) as a pattern, we can write
\begin{equation}
   Z = LZ \oplus B, \qquad Y = DY \oplus Z, \qquad X = MX \oplus Y.
\label{LDM}
\end{equation}
Then
\begin{equation}
   A^* = M^*D^*L^*.
\label{AMDLstar}
\end{equation}
A triple $(L,D,M)$ consisting of a lower triangular, diagonal, and
upper triangular matrices is called an $LDM$-{\it factorization} of a
matrix $A$ if relations (10) and (11) are satisfied. We note that
in this case, the principal diagonals of $L$ and $M$ are zero.

The modification of the notion of $LDM$-factorization used in matrix
analysis for the equation $AX=B$ is constructed in analogy with the
construct suggested by Carr\'e in [10, 11] for $LU$-factorization.

We stress that the algorithm described below can be applied to matrix
computations over any semiring under the condition that the unary
operation $a\mapsto a^*$ is applicable every time it is encountered
in the computational process. Indeed, when constructing the
algorithm, we use only the basic semiring operations of addition
$\oplus$ and multiplication $\odot$ and the properties of
associativity, commutativity of addition, and distributivity of
multiplication over addition.

If $A$ is a symmetric matrix over a semiring with a commutative
multiplication, the amount of computations can be halved, since
$M$ and $L$ are mapped into each other under transposition.

We begin with the case of a triangular matrix $A = L$ (or $A = M$).
Then, finding $X$ is reduced to the forward (or back) substitution.

\begin{center}
{\it Forward substitution}
\end{center}

 We are given:
\begin{itemize}
\item $L = \|l^i_j\|^n_{i,j = 1}$, where $l^i_j = \0$ for $i \le j$
(a lower triangular matrix with a zero diagonal);
\item $B = \|b^i\|^n_{i = 1}$.
\end{itemize}

It is required to find the solution $X = \|x^i\|^n_{i = 1}$ to the
equation $X = LX \oplus B$. The program fragment solving this problem is as
follows.

\begin{tabbing}
   \qquad\=\qquad\=\kill
   for $i = 1$ to $n$ do\\*
   \{\> $x^i := b^i$;\\
   \> for $j = 1$ to $i - 1$ do\\*
   \>\> $x^i := x^i \oplus (l^i_j \odot x^j)$;\, \}\\
\end{tabbing}

\begin{center}
{\it Back substitution}
\end{center}

We are given
\begin{itemize}
\item $M = \|m^i_j\|^n_{i,j = 1}$, where $m^i_j = \0$ for $i \ge j$ (an
upper triangular matrix with a zero diagonal);
\item $B = \|b^i\|^n_{i = 1}$.
\end{itemize}

It is required to find the solution $X = \|x^i\|^n_{i = 1}$ to the
equation $X = MX \oplus B$. The program fragment solving this problem
is as follows.

\begin{tabbing}
   \qquad\=\qquad\=\kill
   for $i = n$ to 1 step $-1$ do\\*
   \{\> $x^i :=  b^i$;\\
   \> for $j = n$ to $i + 1$ step $-1$ do\\*
   \>\> $x^i :=  x^i \oplus (m^i_j \odot x^i)$;\, \}\\
\end{tabbing}

Both algorithms require $(n^2 - n) / 2$ operations $\oplus$ and $\odot$.

\begin{center}
{\it Closure of a diagonal matrix}
\end{center}

We are given
\begin{itemize}
\item $D = {\rm{diag}}(d_1, \ldots, d_n)$;
\item $B = \|b^i\|^n_{i = 1}$.
\end{itemize}

It is required to find the solution $X = \|x^i\|^n_{i = 1}$ to the
equation $X = DX \oplus B$. The program fragment solving this problem
is as follows.

\begin{tabbing}
   \qquad\=\qquad\=\kill
   for $i = 1$ to $n$ do\\*
   \> $x^i :=  (d_i)^* \odot b^i$;\\
\end{tabbing}

This algorithm requires $n$ operations $*$ and $n$ multiplications $\odot$.

\begin{center}
{\it General case}
\end{center}

We are given

\begin{itemize}
\item $L = \|l^i_j\|^n_{i,j = 1}$, where $l^i_j = \0$ if $i \le j$;
\item $D = {\rm{diag}}(d_1, \ldots, d_n)$;
\item $M = \|m^i_j\|^n_{i,j = 1}$, where $m^i_j = \0$ if $i \ge j$;
\item $B = \|b^i\|^n_{i = 1}$.
\end{itemize}

It is required to find the solution $X = \|x^i\|^n_{i = 1}$ to the
equation $X = AX \oplus B$, where $L$, $D$, and $M$ form the
$LDM$-factorization of $A$. The program fragment solving this problem
is as follows.

\begin{tabbing}
        {\sc {FORWARD SUBSTITUTION}}\\*
   for $i = 1$ to $n$ do\\*
   \{\, $x^i :=  b^i$;\\*
   \, for $j = 1$ to $i - 1$ do\\*
   \,\, $x^i :=  x^i \oplus (l^i_j \odot x^j)$;\, \}\\
	\sc{CLOSURE OF A DIAGONAL MATRIX}\\*
   for $i = 1$ to $n$ do\\*
   \, $x^i :=  (d_i)^* \odot b^i$;\\
	\sc{BACK SUBSTITUTION}\\*
   for $i = n$ to 1 step $-1$ do\\*
   \{\, for $j = n$ to $i + 1$ step $-1$ do\\*
   \,\, $x^i :=  x^i \oplus (m^i_j \odot x^j)$;\, \}\\
\end{tabbing}

Note that $x^i$ is not initialized in the course of the back substitution.
The algorithm requires $n^2 - n$ operations $\oplus$, $n^2$ operations
$\odot$, and $n$ operations~$*$.

\begin{center}
{\it LDM-factorization}
\end{center}

We are given
\begin{itemize}
\item $A = \|a^i_j\|^n_{i,j = 1}$.
\end{itemize}

It is required to find the $LDM$-factorization of $A$:
$L = \|l^i_j\|^n_{i,j = 1}$, $D ={\rm{diag}}(d_1, \ldots, d_n)$, and
$M = \|m^i_j\|^n_{i,j = 1}$, where $l^i_j = \0$ if $i \le j$, and
$m^i_j = \0$ if $i \ge j$.

The program uses the following internal variables:
\begin{itemize}
\item $C = \|c^i_j\|^n_{i,j = 1}$;
\item $V = \|v^i\|^n_{i = 1}$;
\item $d$.
\end{itemize}

\begin{tabbing}
   \qquad\=\qquad\=\qquad\=\kill
   \sc{INITIALISATION}\\*
	for $i = 1$ to $n$ do\\*
	\> for $j = 1$ to $n$ do\\*
	\>\> $c^i_j = a^i_j$;\\
	\sc{MAIN LOOP}\\*
	for $j = 1$ to $n$ do\\*
	\{\> for $i = 1$ to $j$ do\\*
	\>\> $v^i :=  a^i_j$;\\
	\> for $k = 1$ to $j - 1$ do\\*
	\>\> for $i = k + 1$ to $j$ do\\*
	\>\>\> $v^i :=  v^i \oplus (a^i_k \odot v^k)$;\\
	\> for $i = 1$ to $j - 1$ do\\*
	\>\> $a^i_j :=  (a^i_i)^* \odot v^i$;\\
	\> $a^j_j :=  v^j$;\\
	\> for $k = 1$ to $j - 1$ do\\*
	\>\> for $i = j + 1$ to $n$ do\\*
	\>\>\> $a^i_j :=  a^i_j \oplus (a^i_k \odot v^k)$;\\
	\> $d = (v^j)^*$;\\
	\> for $i = j + 1$ to $n$ do\\*
	\>\> $a^i_j :=  a^i_j \odot d$;\, \}\\
\end{tabbing}

This algorithm requires $(2n^3 - 3n^2 + n) /6$ operations $\oplus$, $(2n^3 +
3n^2 -5n) / 6$ operations $\odot$, and $n(n + 1) / 2$ operations $*$.
After its completion, the matrices $L$, $D$, and $M$ are contained,
respectively, in the lower triangle, on the diagonal, and in the upper
triangle of the matrix $C$. In the case when $A$ is symmetric about the
principal diagonal and the semiring over which the matrix is defined
is commutative, the algorithm can be modified in such a way that the
number of operations is reduced approximately by a factor of two.

Other examples can be found in [10, 11, 22--25, 35, 36, 66, 67].

Note that to compute the matrices $A^*$ and $A^*B$ it is convenient to solve the Bellman equation (2).

Some other interesting and important problems of linear algebra over
semirings are examined, e.g., in [8, 9, 12, 18, 20, 22--25, 27, 57, 58, 59, 60, 74--77].

\section{The idempotent correspondence principle}

There is a nontrivial analogy between mathematics of semirings and
quantum mechanics. For example, the field of real numbers can be
treated as a ``quantum object'' with respect to idempotent semirings.
So idempotent semirings can be treated as ``classical'' or
``semi-classical'' objects with respect to the field of real
numbers.

Let $\rset$ be the field of real numbers and $\rset_+$ the subset
of all non-negative numbers. Consider the following change of variables:
$$
u \mapsto w = h \ln u,
$$
where $u \in \rset_+ \setminus \{0\}$, $h > 0$; thus $u = e^{w/h}$,
$w \in
\rset$. Denote by $\0$ the additional element $-\infty$ and by $S$
the extended real line $\rset \cup \{\0\}$. The above change of
variables has a
natural extension $D_h$ to the whole $S$ by $D_h(0) = \0$; also, we
denote $D_h(1) = 0 = \1$.

Denote by $S_h$ the set $S$ equipped with the two operations $\oplus_h$
(generalized addition) and $\odot_h$ (generalized multiplication)
such that
$D_h$ is a homomorphism of $\{\rset_+, +, \cdot\}$ to $\{S, \oplus_h,
\odot_h\}$. This means that $D_h(u_1 + u_2) = D_h(u_1) \oplus_h D_h(u_2)$
and $D_h(u_1 \cdot u_2) = D_h(u_1) \odot_h D_h(u_2)$, i.e., $w_1 \odot_h
w_2 = w_1 + w_2$ and $w_1 \oplus_h w_2 = h \ln (e^{w_1/h} + e^{w_2/h})$.
It is easy to prove that $w_1 \oplus_h w_2 \to \max\{w_1, w_2\}$
as $h \to 0$.

Denote by $\rmax$ the set $S = \rset \cup \{\0\}$ equipped
with operations
$\oplus = \max$ and $\odot = +$, where ${\0} = -\infty$, ${\1} = 0$
as above.
Algebraic structures in $\rset_+$ and $S_h$ are isomorphic; therefore
$\rmax$ is a result of a deformation of the structure in $\rset_+$.

We stress the obvious analogy with the quantization procedure,
where $h$ is
the analog of the Planck constant. In these terms, $\rset_+$
(or $\rset$)
plays the part of a ``quantum object'' while $\rmax$ acts as a
``classical'' or ``semi-classical'' object that arises as the result
of a {\it dequantization} of this quantum object.

Likewise, denote by $\rmin$ the set $\rset \cup \{\0\}$ equipped with
operations $\oplus = \min$ and $\odot = +$, where ${\0} = +\infty$ and
${\1} = 0$. Clearly, the corresponding dequantization procedure is
generated by the change of variables $u \mapsto w = -h \ln u$.

Consider also the set $\rset \cup \{\0, \1\}$, where ${\0} = -\infty$,
${\1} =+\infty$, together with the operations $\oplus = \max$
and $\odot=\min$.
Obviously, it can be obtained as a result of a ``second dequantization''
of $\rset$ or $\rset_+$.

There is a natural transition from the field of real numbers or complex numbers to the semiring $\rmax$ (or $\rmin$). This is a composition of the mapping $x\mapsto |x|$ and the deformation described above.

In general an \emph{idempotent dequantization} is a transition from a basic field to an idempotent semiring in mathematical concepts, constructions and results, see [38--40] for details.

For example, the basic object of the traditional calculus is a {\it function}
defined on some set $X$ and taking its values in the field $\rset$
(or $\textbf{C}$); its idempotent analog is a map $X \to S$, where $X$ is
some set and $S =\rmin$, $\rmax$, or another idempotent semiring. Let us
show that redefinition of basic constructions of traditional calculus in
terms of idempotent mathematics can yield interesting and nontrivial
results (see, e.g., [32, 37--43, 45, 46, 48--53, 71, 72], for details,
additional examples and generalizations).

\noindent\textbf{Example~6. Integration and measures.} 
To define an idempotent
analog of the Riemann integral, consider a Riemann sum for a function
$\varphi(x)$, $x \in X = [a,b]$, and substitute semiring
operations $\oplus$
and $\odot$ for operations $+$ and $\cdot$ (usual addition and
multiplication) in its expression (for the sake of being definite,
consider the semiring $\rmax$):
$$
\sum_{i = 1}^N \varphi(x_i) \cdot \Delta_i \quad\mapsto\quad
\bigoplus_{i = 1}^N \varphi(x_i) \odot \Delta_i
= \max_{i = 1, \ldots, N}\, (\varphi(x_i) + \Delta_i),
$$
where $a = x_0 < x_1 < \cdots < x_N = b$, $\Delta_i = x_i - x_{i - 1}$, $i
= 1,\ldots,N$. As $\max_i \Delta_i \to 0$, the integral sum tends to
$$
 \int_X^\oplus \varphi(x)\, dx = \sup_{x \in X} \varphi(x)
$$
for any function $\varphi$:~$X \to \rmax$ that is bounded. In general,
for any set $X$ the set function
$$
m_\varphi(B) = \sup_{x \in B} \varphi(x), \quad B \subset X,
$$
is called an $\rmax$-{\it measure} on $X$;
since $m_\varphi(\bigcup_\alpha
B_\alpha) = \sup_\alpha m_\varphi(B_\alpha)$, this measure is completely
additive. An idempotent integral with respect to this measure
is defined as
$$
 \int_X^\oplus \psi(x)\, dm_\varphi
 = \int_X^\oplus \psi(x) \odot \varphi(x)\, dx
 = \sup_{x \in X}\, (\psi(x) + \varphi(x)).
$$

Using the standard partial order it is possible to generalize these
definitions for the case of arbitrary idempotent semirings.

\noindent\textbf{Example~7. Fourier--Legendre transform.}
Consider the topological
group $G = \rset^n$. The usual Fourier--Laplace transform is defined as
$$
\varphi(x) \mapsto \widetilde\varphi(\xi)
 = \int_G e^{i\xi \cdot x} \varphi(x)\, dx,
$$
where $\exp(i\xi \cdot x)$ is a {\it character} of the group $G$, i.e.,
a solution of the following functional equation:
$$
   f(x + y) = f(x)f(y).
$$

The idempotent analog of this equation is
$$
 f(x + y) = f(x) \odot f(y) = f(x) + f(y).
$$
Hence natural ``idempotent characters'' of the group $G$ are linear
functions of the
form $x \mapsto \xi \cdot x = \xi_1 x_1 + \cdots + \xi_n x_n$. Thus
the Fourier--Laplace transform turns into
$$
 \varphi(x) \mapsto \widetilde\varphi(\xi)
 = \int_G^\oplus \xi \cdot x \odot \varphi(x)\, dx
 = \sup_{x \in G}\, (\xi \cdot x + \varphi(x)).
$$
This is the well-known Legendre (or Fenchel) transform. Examples related to an important version of matrix algebra are discussed in section 4 above.

These examples suggest the following formulation of the idempotent
correspondence principle [39, 40]:
\begin{quote}
{\it There exists a heuristic correspondence between interesting, useful,
and important constructions and results over the field of real (or
complex) numbers and similar constructions and results over idempotent
semirings in the spirit of N.~Bohr's correspondence principle in
quantum mechanics.}
\end{quote}

So idempotent mathematics can be treated as a ``classical shadow (or
counterpart)'' of the traditional Mathematics over fields.

A systematic application of this correspondence principle leads to a variety of theoretical and applied results, see, e.g. [37--43, 45, 46, 55, 71, 72]. Relations between idempotent and traditional mathematics are presented in Fig 1. Relations to quantum physics are discussed in detail, e.g., in [38].

\begin{figure}
\begin{center}
\includegraphics[width=100 mm,keepaspectratio]{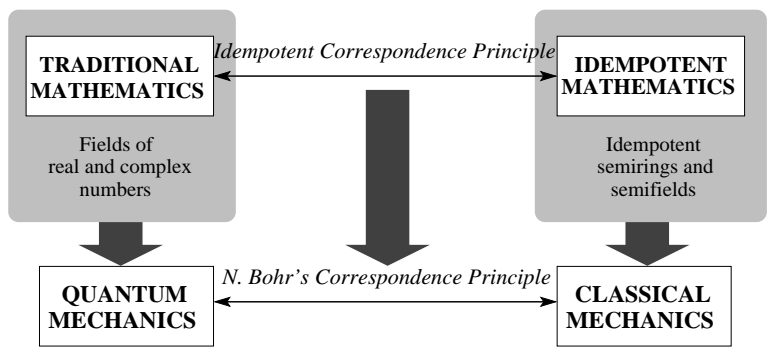}
\caption{Relations between idempotent and traditional mathematics.}
\end{center}
\end{figure}

\section{The superposition principle and parallel computing}

In quantum mechanics the superposition principle means that the
Schr\"odi\-n\-ger equation (which is basic for the theory) is linear.
Similarly in idempotent mathematics the idempotent superposition
principle means that some important and basic problems and equations
(e.g., optimization problems, the Bellman equation and its versions
and generalizations, the Hamilton--Jacobi equation) nonlinear in the
usual sense can be  treated as linear over appropriate idempotent
semirings (in a general form this superposition principle was formulated by V.~P.~Maslov), see [48--53].

Linearity of the Hamilton-Jacobi equation over $\rset_{\min}$
(and $\rset_{\max}$) can be deduced from the usual linearity (over
$\textbf{C}$) of the corresponding Schr\"odinger equation by means of the
dequantization procedure described above (in Section 4). In this case
the parameter $h$ of this dequantization coincides with $i\hbar$ ,
where $\hbar$ is the Planck constant; so in this case $\hbar$ must take
imaginary values (because $h>0$; see [38] for details).
Of course, this is closely related to variational principles of
mechanics.

The situation is similar for the differential Bellman equation (see, e.g.,
[32]) and the discrete version of the Bellman equations, see section 4 above.

It is well known that linear problems and equations are especially convenient for parallelization, see, e.g., [73].
Standard methods (including the so-called block methods) constructed in the framework of the traditional mathematics can be extended to universal algorithms over semirings (the correspondence principle!). For example, formula (1) discussed in section 4.3 leads to a simple block method for parallelization of the closure operations. Other standard methods of linear algebra [73] can be used in a similar way.

\section{The correspondence principle for computations}

Of course, the idempotent correspondence principle is valid for
algorithms as well as for their software and hardware implementations
[39, 40, 42, 43]. Thus:

\begin{quote}
{\it If we have an important and interesting numerical algorithm, then
there is a good chance that its semiring analogs are important and
interesting as well.}
\end{quote}

In particular, according to the superposition principle,
analogs of linear
algebra algorithms are especially important. Note that
numerical algorithms
for standard infinite-dimensional linear problems over idempotent
semirings (i.e., for
problems related to idempotent integration, integral operators and
transformations, the Hamilton-Jacobi and generalized Bellman equations)
deal with the corresponding finite-dimensional (or finite) ``linear
approximations''. Nonlinear algorithms often can be approximated by linear
ones. Thus the idempotent linear algebra is a basis for the idempotent
numerical analysis.

Moreover, it is well-known that linear algebra algorithms easily lend themselves to parallel computation; their idempotent analogs admit
parallelization as
well. Thus we obtain a systematic way of applying parallel computing to
optimization problems.

Basic algorithms of linear algebra (such as inner product of two vectors,
matrix addition and multiplication, etc.) often do not depend on
concrete semirings, as well as on the nature of domains containing the
elements of vectors and matrices. Algorithms to construct the closure
$A^*={\1}\oplus A\oplus A^2\oplus\cdots\oplus A^n\oplus\cdots=
\bigoplus^{\infty}_{n=1} A^n$ of an idempotent matrix $A$ can be derived
from standard methods for calculating $({\1} -A)^{-1}$. For the
Gauss--Jordan elimination method (via LU-decomposition) this trick was used in [66],
and the corresponding algorithm is universal and can be applied both to
the Bellman equation and to computing the inverse of a real (or complex)
matrix $({\1} - A)$. Computation of $A^{-1}$ can be derived
from this universal
algorithm with some obvious cosmetic transformations.

Thus it seems reasonable to develop universal algorithms that can deal
equally well with initial data of different domains sharing the same
basic structure [39, 40, 43].

\section{The correspondence principle for hardware design}

A systematic application of the correspondence principle to computer
calculations leads to a unifying approach to software and hardware
design.

The most important and standard numerical algorithms have many hardware
realizations in the form of technical devices or special processors.
{\it These devices often can be used as prototypes for new hardware
units generated by substitution of the usual arithmetic operations
for its semiring analogs and by addition tools for performing neutral
elements $\0$ and} $\1$ (the latter usually is not difficult). Of course,
the case of numerical semirings consisting of real numbers (maybe except
neutral elements)  and semirings of numerical intervals is the most simple and natural [38--43, 45, 46, 69].
Note that for semifields (including $\rmax$ and $\rmin$)
the operation of division is also defined.

Good and efficient technical ideas and decisions can be transposed
from prototypes into new hardware units. Thus the correspondence
principle generated a regular heuristic method for hardware design.
Note that to get a patent it is necessary to present the so-called
`invention formula', that is to indicate a prototype for the suggested
device and the difference between these devices.

Consider (as a typical example) the most popular and important algorithm
of computing the scalar product of two vectors:
\begin{equation}
(x,y)=x_1y_1+x_2y_2+\cdots + x_ny_n.
\end{equation}
The universal version of (12) for any semiring $A$ is obvious:
\begin{equation}
(x,y)=(x_1\odot y_1)\oplus(x_2\odot y_2)\oplus\cdots\oplus
(x_n\odot y_n).
\end{equation}
In the case $A=\rmax$ this formula turns into the following one:
\begin{equation}
(x,y)=\max\{ x_1+y_1,x_2+y_2, \cdots, x_n+y_n\}.
\end{equation}

This calculation is standard for many optimization algorithms, so
it is useful to construct a hardware unit for computing (14). There
are many different devices (and patents) for computing (12) and every
such device can be used as a prototype to construct a new device for
computing (14) and even (13). Many processors for matrix multiplication
and for other algorithms of linear algebra are based on computing
scalar products and on the corresponding ``elementary'' devices
respectively, etc.

There are some methods to make these new devices more universal than
their prototypes. There is a modest collection of possible operations
for standard numerical semirings: max, min, and the usual arithmetic
operations. So, it is easy to construct programmable hardware
processors with variable basic operations. Using modern technologies
it is possible to construct cheap special-purpose multi-processor
chips implementing examined algorithms. The so-called
systolic processors are
especially convenient for this purpose. A systolic array is a
`homogeneous' computing medium consisting of elementary
processors, where the general scheme and processor connections
are simple and regular. Every elementary processor pumps data in and
out performing elementary operations in a such way that the
corresponding data flow is kept up in the computing medium; there
is an analogy with the blood circulation and this is a reason for the
term ``systolic'', see e.g., [39, 40, 43, 52, 65--67].

Concrete systolic processors for the general algebraic path problem are
presented in [65--67]. In particular, there is a systolic array of
$n(n+1)$ elementary processors which performs computations of the Gauss--Jordan
elimination algorithm and can solve the algebraic path problem within $5n-2$
time steps. Of course, hardware implementations for important and popular basic
algorithms increase the speed of data processing.

The so-called GPGPU (General-Purpose computing on Graphics Processing Units) technique is another important field for applications of the correspondence principle. The thing is that graphic processing units (hidden in modern laptop and desktop computers) are potentially powerful processors for solving numerical problems. The recent tremendous progress in graphical processing hardware and software leads to new ``open'' programmable parallel computational devices (special processors), see, e.g., [78--80]. These devices are going to be standard for coming PC (personal computers) generations. Initially used for graphical processing only (at that time they were called GPU), today they are used for various fields, including audio and video processing, computer simulation, and encryption.   But this list can be considerably enlarged following the correspondence principle: the basic operations would be used as parameters. Using the technique described in this paper (see also our references), standard linear algebra algorithms can be used for solving different problems in different areas. In fact, the hardware supports all operations needed for the most important idempotent semirings: plus, times, min, max. The most popular linear algebra packages [ATLAS (Automatically Tuned Linear Algebra Software), LAPACK, PLASMA (Parallel Linear Algebra for Scalable Multicore Architectures)] can already use GPGPU, see [81--83]. We propose to make these tools more powerful by using parameterized algorithms.

Linear algebra over the most important numerical semirings generates solutions for many concrete problems in different areas, see, e.g.,
sections 4.4 and 4.5 above and references indicated in these sections.

Note that to be consistent with operations we have to redefine zero (0) and unit (1) elements (see above); comparison operations must be also redefined as it is described in section~4.1 ``Basic definitions.'' Once the operations are redefined, then the most of basic linear algebra algorithms, including back and forward substitution, Gauss elimination method, Jordan elimination method and others could be rewritten for new domains and data structures. Combined with the power of the new parallel hardware this approach could change PC from entertainment devices to power full instruments.

\section{The correspondence principle for software design}

Software implementations for universal semiring algorithms are not
as efficient as hardware ones (with respect to the computation speed)
but they are much more flexible. Program modules can deal with abstract (and
variable) operations and data types. Concrete values for these
operations and data types can be defined by the corresponding
input data. In this case concrete operations and data types are generated
by means of additional program modules. For programs written in
this manner it is convenient to use special techniques of the
so-called object oriented (and functional) design, see, e.g.,
[47, 62, 70]. Fortunately, powerful tools supporting the
object-oriented software design have recently appeared including compilers
for real and convenient programming languages (e.g. $C^{++}$ and Java) and modern computer algebra systems.

Recently, this type of programming technique has been dubbed
generic programming (see, e.g., [6, 62]). To help automate the
generic programming, the so-called Standard Template Library (STL)
was developed in the framework of $C^{++}$ [62, 70]. However,
high-level tools, such as STL, possess both obvious advantages
and some disadvantages and must be used with caution.

It seems that it is natural to obtain an implementation of the correspondence
principle approach to scientific calculations in the form of a
powerful software system based on a collection of universal
algorithms. This approach ensures a working time reduction for
programmers and users because of the software unification.
The arbitrary necessary accuracy and safety of numeric calculations can be ensured
as well.

The system has to contain several levels (including programmer and
user levels) and many modules.

Roughly speaking, it must be divided into three parts. The first part
contains modules that implement domain
modules (finite representations of
basic mathematical objects). The second part implements universal
(invariant) calculation methods. The third part contains modules
implementing model dependent algorithms. These modules may be
used in user programs written in $C^{++}$, Java, Maple, Mathlab etc.

The system has to contain the following modules:

\medskip
\begin{itemize}
\item[---] Domain modules:
\begin{itemize}
\item infinite precision integers;
\item rational numbers;
\item finite precision rational numbers (see [44]);
\item finite precision complex rational numbers;
\item fixed- and floating-slash rational numbers;
\item complex rational numbers;
\item arbitrary precision floating-point real numbers;
\item arbitrary precision complex numbers;
\item $p$-adic numbers;
\item interval numbers;
\item ring of polynomials over different rings;
\item idempotent semirings;
\item interval idempotent semirings;
\item and others.
\end{itemize}
\item[---] Algorithms:
\begin{itemize}
\item linear algebra;
\item numerical integration;
\item roots of polynomials;
\item spline interpolations and approximations;
\item rational and polynomial interpolations and approximations;
\item special functions calculation;
\item differential equations;
\item optimization and optimal control;
\item idempotent functional analysis;
\item and others.
\end{itemize}
\end{itemize}

This software system may be especially useful for designers
of algorithms, software engineers, students and mathematicians.

Note that there are some software systems oriented to calculations with idempotent semirings like $\rmax$; see, e.g., [64]. However these systems do not support universal algorithms.

\section{Interval analysis in idempotent mathematics}

Traditional interval analysis is a nontrivial and popular mathematical area, see, e.g., [3, 20, 34, 54, 56, 59]. An ``idempotent'' version of interval analysis (and moreover interval analysis over positive semirings) appeared in [45, 46, 69]. Later appeared rather many publications on the subject, see, e.g., [12, 20, 27, 57, 58, 77]. Interval analysis over the positive semiring $\textbf{R}_+$ was discussed in [8].

Let a set~$S$ be partially ordered by a relation $\cle$.
A \emph{closed interval} in~$S$ is a subset of the form $\x = [\lx, \ux] =
\{\, x \in S \mid \lx \cle x \cle \ux\, \}$, where the elements $\lx \cle
\ux$ are called \emph{lower} and \emph{upper bounds} of the interval $\x$.
The order~$\cle$ induces a partial ordering on the set of all closed
intervals in~$S$: $\x \cle \y$ iff $\lx \cle \ly$ and $\ux \cle \uy$.

A \emph{weak interval extension} $I(S)$ of a positive semiring~$S$ is the
set of all closed intervals in~$S$ endowed with operations $\oplus$
and~$\odot$ defined as ${\x \oplus \y} = [{\lx \oplus \ly}, {\ux \oplus
\uy}]$, ${\x \odot \y} = [{\lx \odot \ly}, {\ux \odot \uy}]$ and a partial
order induced by the order in $S$. The closure operation in $I(S)$ is
defined by $\x^* = [\lx^*, \ux^*]$. There are some other interval extensions (including the so-called strong interval extension [46]) but the weak extension is more convenient.

The extension $I(S)$ is positive; $I(S)$ is idempotent if $S$ is an idempotent semiring.
A universal algorithm over $S$ can be applied to $I(S)$ and we shall get an interval version of the initial algorithm.
Usually both the versions have the same complexity. For the discrete stationary Bellman equation and the corresponding optimization problems on graphs interval analysis was examined in [45, 46] in details. Other problems of idempotent linear algebra were examined in [12, 27, 57, 58, 77].

Idempotent mathematics appears to be remarkably simpler than its
traditional analog. For example, in traditional interval arithmetic,
multiplication of intervals is not distributive with respect to addition of
intervals, whereas in idempotent interval arithmetic this distributivity is
preserved. Moreover, in traditional interval analysis the set of all
square interval matrices of a given order does not form even a semigroup
with respect to matrix multiplication: this operation is not associative
since distributivity is lost in the traditional interval arithmetic. On the
contrary, in the idempotent (and positive) case associativity is preserved. Finally, in
traditional interval analysis some problems of linear algebra, such as
solution of a linear system of interval equations, can be very difficult
(generally speaking, they are $NP$-hard, see
[14, 20, 33, 34] and references therein). It was noticed  in [45, 46] that in the idempotent case solving an interval linear system
requires a polynomial number of operations (similarly to the usual Gauss
elimination algorithm).  Two properties that make the idempotent interval
arithmetic so simple are monotonicity of arithmetic operations and
positivity of all elements of an idempotent semiring.

Usually interval estimates in idempotent mathematics are exact. In the traditional theory such estimates tend to be overly pessimistic.

\section*{Acknowledgments}

This work is partially supported by the RFBR grant 08-01-00601.
The authors are grateful to A.~G.~Kushner for his kind help.

\end{document}